\documentclass[12pt,reqno]{amsart}

\newcommand\version{September 28, 2021}

\usepackage{amsmath,amsfonts,amsthm,amssymb,amsxtra}
\usepackage{bbm} % to get \1
\usepackage{hyperref}	% to get hyperlinks in equations and references
\newtheorem{theorem}{Theorem}%[section]

\newtheorem{lemma}[theorem]{Lemma}

\theoremstyle{definition}

\theoremstyle{remark}

\newtheorem{remarks}[theorem]{Remarks}

\newcommand{\1}{\mathbbm{1}}

\newcommand{\C}{\mathbb{C}}

\newcommand{\const}{\mathrm{const}\ }

\renewcommand{\epsilon}{\varepsilon}

\newcommand{\N}{\mathbb{N}}

\renewcommand{\phi}{\varphi}
\newcommand{\R}{\mathbb{R}}

\newcommand{\Sph}{\mathbb{S}}

\DeclareMathOperator{\dist}{dist}

\DeclareMathOperator{\dom}{dom}

\DeclareMathOperator{\supp}{supp}

\DeclareMathOperator{\Tr}{Tr}

\begin{document}

\title[Lieb--Thirring inequalities --- \version]{Lieb--Thirring inequalities and other functional inequalities for orthonormal systems}

\author{Rupert L. Frank}
\address[Rupert L. Frank]{Mathe\-matisches Institut, Ludwig-Maximilans Universit\"at M\"unchen, The\-resienstr.~39, 80333 M\"unchen, Germany, and Munich Center for Quantum Science and Technology, Schel\-ling\-str.~4, 80799 M\"unchen, Germany, and Mathematics 253-37, Caltech, Pasa\-de\-na, CA 91125, USA}
\email{r.frank@lmu.de}

\renewcommand{\thefootnote}{${}$} \footnotetext{\copyright\, 2021 by the author. This paper may be reproduced, in its entirety, for non-commercial purposes.}

\begin{abstract}
	We review recent results on functional inequalities for systems of orthonormal functions. The key finding is that for various operators the orthonormality leads to a gain over a simple application of the triangle inequality. The operators under consideration are either related to Sobolev type inequalities or to Fourier restriction type inequalities.
\end{abstract}

\maketitle

\section{Introduction}

For more than four decades, Lieb--Thirring inequalities have played an important role in various areas of mathematical physics and analysis. The progress that has been made towards the conjectures in the area and many extensions and generalizations of the original inequalities have been reviewed in the surveys \cite{Li89,BlaStu-96,LaWe1,Hu,La,Fr}, the textbooks \cite{LiLo,LiSe}, as well as in the forthcoming book \cite{FrLaWe}. In order to avoid too large an overlap with these existing works, the present contribution, which was invited by the organizers of the International Congress of Mathematicians 2022, to whom the author is most grateful, focuses only on one single aspect of these inequalities. Namely, we will consider Lieb--Thirring inequalities from the point of view of Sobolev-type inequalities for systems of orthonormal functions, and we discuss recent extensions, in particular, to the Strichartz and Stein--Tomas inequalities from harmonic analysis. We will also briefly comment on selected applications of these newly obtained bounds.

\subsection{The general set-up}\label{sec:setting}

Let $\mathcal H$ be a (typically complex) Hilbert space with norm denoted by $\|\cdot\|$ and let $X$ be a measure space, with measure denoted simply by $dx$ and with corresponding Lebesgue spaces $L^q(X)$. Assume that $T$ is a bounded linear operator from $\mathcal H$ to $L^q(X)$ for some $q>2$. That is, for all $f\in\mathcal H$,
\begin{equation}
	\label{eq:tbound}
	\int_X |Tf|^q\,dx \lesssim \| f\|^q \,. 
\end{equation}
As a consequence, if $f_1,\ldots,f_N$ are normalized in $\mathcal H$, then
\begin{equation*}
	%\label{eq:tboundtriangle}
	\int_X \left( \sum_{n=1}^N |Tf_n|^2 \right)^{q/2}dx \lesssim N^{q/2}.
\end{equation*}
This is a consequence of \eqref{eq:tbound} and the triangle inequality in $L^{q/2}$. The power $N^{q/2}$ is best possible, as can be seen by taking all $f_n$ to be equal.

The question that interests us here is whether for a given operator $T$ there is a power
$$
\alpha<q/2
$$
such that for all $N$ and all $f_1,\ldots,f_N\in\mathcal H$ satisfying the \emph{orthonormality constraint}
$$
(f_n,f_m) = \delta_{n,m}
\qquad\text{for all}\ 1\leq n,m\leq N
$$
one has
\begin{equation}
	\label{eq:tboundlt}
	\int_X \left( \sum_{n=1}^N |Tf_n|^2 \right)^{q/2}dx \lesssim N^\alpha.
\end{equation}
As explained, for instance, in \cite{LiSe,LeLiSe,Fr}, such bounds, if true, have important consequences in the mathematical physics of large fermionic quantum systems, density functional theory and the theory of nonlinear evolution equations. Their study is also interesting from a purely analytical point of view and reveals aspects of the underlying operator $T$ which go beyond its boundedness.

At the moment there is no general principle that determines the exponent $\alpha$ directly from the operator $T$. Rather, inequalities of the form \eqref{eq:tboundlt}, if true, have been verified on a case by case basis. Most of the existing results concern the case where $T$ is (at least approximately) translation invariant. Finding a regime of orthonormal functions $f_1,\ldots,f_N$ with $N\to\infty$ where the power $\alpha$ in the bound \eqref{eq:tboundlt} is saturated relies often on techniques of semiclassical and microlocal analysis.

%%%%%%%%%%%%%%%%

\subsection{Example: The HLS inequality}

The above principle is most clearly illustrated on the example of the Hardy--Littlewood--Sobolev (HLS) inequality, also know as the weak Young inequality or the theorem of fractional integration; see, e.g., \cite[Theorem 4.3]{LiLo}. This inequality states that for $0<s<d/2$ the operator of convolution with $|x|^{-d+s}$ is bounded from $L^2(\R^d)$ to $L^q(\R^d)$ with $q=2d/(d-2s)$.

Its extension to systems of orthonormal functions is due to Lieb \cite{Li-83a} and reads as follows.

\begin{theorem}\label{liebhls}
	Let $0<s<d/2$. Then, if $f_1,\ldots,f_N$ are orthonormal in $L^2(\R^d)$,
	$$
	\int_{\R^d} \left( \sum_{n=1}^N \left| |x|^{-d+s}* f_n \right|^2 \right)^{d/(d-2s)}dx \lesssim N \,.
	$$
\end{theorem}

\begin{remarks}\label{rem:liebhls}
	(a) The power $1$ of $N$ on the right side is best possible.\\
	(b) The bound is equivalent (in a certain weak sense) to the Cwikel--Lieb--Rozenblum (CLR) bound
	$$
	N((-\Delta)^s +V) \lesssim \int_{\R^d} V_-^{d/(2s)}dx
	$$
	on the number $N((-\Delta)^s +V)$ of negative eigenvalues of the generalized Schr\"odinger operator $(-\Delta)^s+V$ in $L^2(\R^d)$. Here $V(x)_-=\max\{-V(x),0\}$. The meaning of `equivalent' will be explained in the next subsection. It is a `weak' form of equivalence, because this argument does \emph{not} mean that the sharp constant in Theorem \ref{liebhls} is in one-to-one correspondence with the sharp constant in the CLR bound. This is in contrast to a form of duality that we will encounter later.\\
	(c) The proof of Theorem \ref{liebhls} in \cite{Li-83a} proceeds by reducing it to Cwikel's proof of the CLR inequality \cite{Cw}. Alternative, direct proofs of Theorem \ref{liebhls} were given in \cite{Ru0,Fr2}. We present a different, unpublished proof in Subsection \ref{sec:liebhls} below.\\
	(d) Just like the HLS inequality, the bound in Theorem \ref{liebhls} is conformally invariant. This leads to a natural conjecture for its optimal constant \cite{Fr2}.
\end{remarks}

%%%%%%%%%%%%%

\subsection{The duality argument}\label{sec:duality}

Let us return to the general setting described in Subsection \ref{sec:setting} and consider a bounded operator $T:\mathcal H\to L^q(X)$ for some $q>2$. By H\"older's inequality, this boundedness is equivalent to having, for any $W\in L^{2q/(q-2)}(X)$ and any $f\in\mathcal H$,
$$
\int_X |W|^2 |T f|^2\,dx \lesssim \| W \|_{2q/(q-2)}^2 \| f\|^2 \,,
$$
which, in turn, is equivalent to the boundedness of the operator $WT$ from $\mathcal H$ to $L^2(X)$ with norm
$$
\| W T \| \lesssim \| W \|_{2q/(q-2)} \,.
$$
Here, as usual, we do not distinguish in the notation between the function $W$ and the operator of multiplication by $W$. Moreover, $\|\cdot\|$ on the left side denotes the operator norm.

\medskip

Let us now reformulate the desired inequality \eqref{eq:tboundlt} in terms of the operator $WT$. We assume that $\alpha<q/2$. Again by H\"older's inequality, we see that \eqref{eq:tboundlt} is equivalent to
\begin{equation}
	\label{eq:duality}
	\sum_{n=1}^N \int_X |W|^2 |Tf_n|^2\,dx \lesssim N^{2\alpha/q} \|W\|_{2q/(q-2)}^2 \,.
\end{equation}

\medskip

In order to state this previous inequality succinctly, we recall the notion of Schatten spaces. Background can be found, for instance, in \cite{GoKr,SiTraceIdeals}. For a compact operator $K$ between two Hilbert spaces, we denote by $(s_n(K))_{n\in\N}$ the sequence of its singular values, that is, the square roots of the eigenvalues of the operator $K^*K$ in nonincreasing order, repeated according to multiplicities. Then, by definition, for any $0<r<\infty$, the Schatten class $\mathcal S^r$ consists of all compact operators $K$ with $s_\cdot(K)\in\ell^r$. This is a normed linear space with respect to
$$
\| K \|_r := \left( \sum_{n\in\N} s_n(K)^r \right)^{1/r}.
$$
Also, we will need the weak variant of this space, $\mathcal S^r_{\mathrm{weak}}$, consisting of all compact $K$ with $s_\cdot(K)\in\ell^r_{\mathrm{weak}}$. For $2<r<\infty$, the following norm will appear naturally in our analysis,
$$
\| K \|_{r, \mathrm{w}} := \sup_{n\in\N} N^{-1/2+1/r} \left( \sum_{n=1}^N s_n(K)^2 \right)^{1/2}.
$$
It follows from the variational principle for sums of eigenvalues that
$$
\sum_{n=1}^N s_n(K)^2 = \sup\left\{ \sum_{n=1}^N \|Kf_n\|^2 :\ f_1,\ldots, f_N\ \text{orthonormal} \right\}.
$$
From this and the triangle inequality in $\R^N$ it follows that $\| \cdot \|_{r, \mathrm{w}}$ defines, indeed, a norm. It is also easy to see that $\| \cdot \|_{r, \mathrm{w}}$ is equivalent to the more standard quasinorm in $\mathcal S^r_{\mathrm{weak}}$ defined by
$$
\| K \|_{r, \mathrm{w}}' := \sup_{n\in\N} n^{1/r} s_n(K) \,.
$$
The constants in this equivalence depend on $r>2$ and their explicit values can be found, for instance, in \cite[Lemma 2.3]{Fr2}, where another expression for $\| K \|_{r, \mathrm{w}}$ is used.

\medskip

Returning to the above setting, we now see that \eqref{eq:duality} is equivalent to the fact that $WT$ belongs to the weak Schatten class $\mathcal S^{2q/(q-2\alpha)}_{\mathrm{weak}}$ with
\begin{equation}
	\label{eq:dualityabstract}
	\| WT \|_{2q/(q-2\alpha),\mathrm{w}} \lesssim \|W\|_{2q/(q-2)} \,.
\end{equation}
To summarize, we have seen that the desired inequality \eqref{eq:tboundlt} is equivalent to a quantitative compactness property of the operator $WT$, expressed in terms of a weak Schatten norm. The exponent $\alpha$ in \eqref{eq:tboundlt} is in one-to-one correspondence with the Schatten exponent. What we have gained through this reformulation is, for instance, that we can use interpolation methods to prove trace ideal properties of the operators $TW$.

\medskip

At this point we can present Lieb's proof of Theorem \ref{liebhls}. Namely, Cwikel's theorem \cite{Cw} says that, for $2<p<\infty$,
\begin{equation}
	\label{eq:cwikel}
	\| a(X) b(-i\nabla) \|_{p,\mathrm{w}} \lesssim \|a\|_p \|b\|_{p,\mathrm{w}} \,.
\end{equation}
Here $a(X)$ denotes the operator of multiplication by a function $a\in L^p(\R^d)$ in position space and $b(-i\nabla)$ denotes the operator of multiplication by a function $b\in L^p_{\mathrm{weak}}(\R^d)$ in momentum space. The operator $T$ relevant for Theorem \ref{liebhls} is convolution with $|x|^{-d+s}$ which corresponds to multiplication by (a constant times) $|\xi|^{-s}$ in Fourier space. The latter function belongs to $L^{d/s}_{\mathrm{weak}}(\R^d)$. Thus, \eqref{eq:cwikel} implies \eqref{eq:dualityabstract} with $\alpha=1$ and $q=2d/(d-2s)$, as claimed.

\medskip

The proof of Cwikel's theorem in \cite{Fr2} goes in some sense the other way around. Namely, first Theorem \ref{liebhls} (or rather a slight generalization of it) is established, using the method of \cite{Ru0}, and then the above duality argument is used to deduced \eqref{eq:cwikel}.

\medskip

We can now also explain the notion of weak equivalence in Remark \ref{rem:liebhls} (b). Namely, by the Birman--Schwinger principle the bounds there for negative eigenvalues of generalized Schr\"odinger operators are the same as bounds on the operator $W(-\Delta)^{-s/2}$ in the quasinorm $\|\cdot\|_{d/s,{\rm weak}}'$, whereas by the above argument the bound in Theorem \ref{liebhls} is the same as a bound on this operator in the norm $\|\cdot\|_{d/s,{\rm weak}}$.

%%%%%%%%%%%%%

\subsection{A generalization}\label{sec:hlsgen}

There is a far reaching generalization of Theorem \ref{liebhls}. Namely, if $X$ is a sigma-finite measure space and $A$ is a nonnegative operator in $L^2(X)$ with heat semigroup satisfying, for some $\nu>2$,
$$
\| \exp(-tA) \|_{L^2\to L^\infty} \lesssim t^{-\nu/4}
\qquad\text{for all}\ t>0 \,,
$$
then for all $u_1,\ldots,u_N\!\in\!\dom A^{1/2}$ satisfying $(A^{1/2}u_n,A^{1/2} u_m)=\delta_{n,m}$ for $1\leq n,m\leq N$,
$$
\int_X \left( \sum_{n=1}^N |u_n|^2 \right)^{\nu/(\nu-2)} dx \lesssim N \,.
$$
This is shown in \cite{Fr2}, improving earlier results in \cite{LeSo,FrLiSe3} that require nonnegativity of the heat kernel.

This more general result reduces to Theorem \ref{liebhls} with $A=(-\Delta)^s$ and $u_n=$ a constant times $(-\Delta)^{-s/2}f_n$. Another application concerns the case where $A$ is the Laplace--Beltrami operator on certain noncompact manifolds. For a compact manifold the above assumption on the semigroup is not satisfied because of the zero eigenvalue, but one can add a positive constant to the Laplace--Beltrami operator.

%%%%%%%%%%%%%

\subsection{Appendix: Proof of Theorem \ref{liebhls}}\label{sec:liebhls}

We present here an unpublished proof of Theorem \ref{liebhls}. It is neither the most elementary one, nor one giving particularly good constants, but we think it is conceptually rather clear and might allow for interesting generalizations. In view of the previous subsection it provides an alternative proof of the CLR inequality and is based on some ideas of Conlon \cite{Co}.

\medskip

By the duality argument in Subsection \ref{sec:duality}, we need to prove \eqref{eq:dualityabstract} with $T$ equal to convolution with $|x|^{-d+s}$, $\alpha=1$ and $q=2d/(d-2s)$. Since the weak Schatten norm of $WT$ equals that of $(WT)^* = T\overline W$, it suffices to consider the latter operator. We have, using $\int |x-z|^{-d+s} |z-y|^{-d+s}\,dz = \const |x-y|^{-d+2s}$ and the Fefferman--de la Llave decomposition \cite{FedlL},
\begin{align*}
	\left\| T \overline W f \right\|^2 & = \const \iint_{\R^d\times\R^d} \frac{\overline{f(x)}\, W(x)\,\overline{W(y)}\,f(y)}{|x-y|^{d-2s}}\,dx\,dy \\
	& = \const \int_{\R^d} da \int_0^\infty \frac{dr}{r^{2d+1-2s}} \iint_{B_r(a)\times B_r(a)} \overline{f(x)}\, W(x)\,\overline{W(y)}\,f(y)\,dx\,dy \,.
\end{align*}
We apply this with $f=f_n$ for some orthonormal $f_n$ in $L^2(\R^d)$ and sum over $n$. For fixed $a$ and $r$, we estimate the double integral over $x$ and $y$ in two different ways. First, since the operator $\gamma$ with kernel $\sum_n \overline{f_n(x)} f_n(y)$ has operator norm one, we have
\begin{align*}
	& \left| \sum_{n=1}^N \iint_{B_r(a)\times B_r(a)} \overline{f_n(x)}\, W(x)\,\overline{W(y)}\,f_n(y)\,dx\,dy \right| = \left| ( \overline W\1_{B_r(a)},\gamma \overline W\1_{B_r(a)}) \right| \\
	& \qquad \leq \int_{B_r(a)} |W(x)|^2\,dx \lesssim r^d (M(|W|^2))(a) \,,
\end{align*}
where $M$ is the maximal function. Second, since $\gamma\geq 0$,
$$
\left| \sum_{n=1}^N \overline{f_n(x)} f_n(y) \right| \leq \sqrt{\rho(x)} \sqrt{\rho(y)} \,,
\qquad \text{where}\ \rho(x):= \sum_{n=1}^N |f_n(x)|^2 \,,
$$
so
\begin{align*}
	& \left| \sum_{n=1}^N \iint_{B_r(a)\times B_r(a)} \overline{f_n(x)}\, W(x)\,\overline{W(y)}\,f_n(y)\,dx\,dy \right| 
	\leq \left( \int_{B_r(a)} \sqrt{\rho(x)}\, |W(x)|\,dx \right)^2 \\
	& \qquad \lesssim r^{2d} (M(|W|\sqrt\rho))(a)^2 \,.
\end{align*}
Inserting this into the above formula, we find
\begin{align*}
	\sum_{n=1}^N \left\| T \overline W f_n \right\|^2 & \lesssim \int_{\R^d} da \int_0^\infty \frac{dr}{r^{2d+1-2s}} \min\left\{ r^d (M(|W|^2))(a), r^{2d} (M(|W|\sqrt\rho))(a)^2 \right\} \\
	& = \const \int_{\R^d} da \, (M(|W|\sqrt\rho))(a)^{2(1-2s/d)} (M(|W|^2))(a)^{2s/d} \\
	& \leq \const \left( \int_{\R^d} (M(|W|\sqrt\rho))(a)^{2d/(d+2s)} da\right)^{1-(2s/d)^2} \\
	& \quad \times \left( \int_{\R^d} (M(|W|^2))(a)^{d/(2s)}\,da \right)^{(2s/d)^2}.
\end{align*}
By the boundedness of the maximal function on $L^p$, $1<p<\infty$, this is bounded by a constant times
\begin{align*}
	& \left( \int_{\R^d} |W|^{2d/(d+2s)} \rho^{d/(d+2s)} da\right)^{1-(2s/d)^2} \left( \int_{\R^d} |W|^{d/s}\,da \right)^{(2s/d)^2} \\
	& \leq \left( \int_{\R^d} |W|^{d/s}\,da \right)^{2s/d} \left( \int_{\R^d} \rho\,da \right)^{1-2s/d}.
\end{align*}
To summarize, we have shown that
$$
\sum_{n=1}^N \left\| T \overline W f_n \right\|^2 \lesssim \|W\|_{d/s}^2 N^{1-2s/d} \,.
$$
If we take the $f_n$ to be the eigenfunctions of $WT^2\overline W$ corresponding to its $N$ largest eigenvalues, the previous inequality becomes
$$
\sum_{n=1}^N s_n(T\overline W)^2 \lesssim \|W\|_{d/s}^2 N^{1-2s/d} \,.
$$
This is the claimed bound on $T\overline W$ in the Schatten space $\mathcal S^{d/s}_\mathrm{weak}$.

%%%%%%%%%%%%%

\section{Sobolev-type inequalities for orthonormal functions}

Before turning to the more recent bounds related to Fourier restriction, in this section we review some classical inequalities for orthonormal functions that are related to Sobolev inequalities. Those include, in particular, the classical Lieb--Thirring inequality in Theorem \ref{lt} below.

\subsection{Bessel-potential bounds}

The bounds in Theorem \ref{liebhls} concern $|x|^{-d+s} *f$, which is a constant multiple of the Riesz potential $(-\Delta)^{-s/2}f$. We present a generalization, due to Lieb \cite{Li-83a}, of these bounds to the Bessel potentials $(-\Delta+m^2)^{-s/2} f$ with $m>0$.

\begin{theorem}\label{liebhlsbessel}
	Let $s>0$ and let 
	$$
	\begin{cases} 2\leq q \leq \infty & \text{if}\ s>d/2 \,,\\
		2\leq q<\infty & \text{if}\ s=d/2 \,,\\
		2\leq q\leq 2d/(d-2s) & \text{if}\ s<d/2 \,.
	\end{cases}
	$$
	Then, if $f_1,\ldots,f_N$ are orthonormal in $L^2(\R^d)$,
	$$
		\left\| \sum_{n=1}^N \left| (-\Delta+m^2)^{-s/2} f_n \right|^2 \right\|_{q/2} \lesssim m^{d-2s-2d/q}  N^{2/q} \,.
	$$
\end{theorem}

The bound for $q=2d/(d-2s)$ if $s<d/2$ follows as before using Cwikel's theorem \eqref{eq:cwikel}. The remaining bounds follow similarly, but using the simpler bound
\begin{equation}
	\label{eq:kss}
	\left\| a(X) b(-i\nabla) \right\|_p \leq (2\pi)^{-d/p} \|a\|_p \|b\|_{p}
\end{equation}
for $2\leq p\leq\infty$. The latter bound is due to Kato, Seiler and Simon (see, e.g., \cite[Theorem 4.1]{SiTraceIdeals}) and can also be inferred from the Lieb--Thirring matrix inequality \cite{LiTh}.

Since \eqref{eq:kss}, in contrast to \eqref{eq:cwikel}, involves a strong instead of a weak Schatten norm on the left side, a generalization of the bound in Theorem \ref{liebhlsbessel} to sums of the form $\sum_n \nu_n \left| (-\Delta+m^2)^{-s/2} f_n \right|^2$ is possible. We discuss this in the next section.

Using bounds due to Solomyak \cite{Sol} (and their natural extension to odd dimensions) it seems plausible that in the case $s=d/2$ there is an endpoint bound in the Orlicz space $\exp L$ in the spirit of a Moser--Trudinger inequality. For instance, the bounds in \cite{FrLa} can be dualized to yield that, if $\Omega\subset\R^2$ is open and of finite measure, then for any $u_1,\ldots, u_N\in H^1_0(\Omega)$ satisfying $\int_\Omega \nabla \overline{u_n}\cdot\nabla u_m\,dx = \delta_{n,m}$ for all $1\leq n,m\leq N$,
$$
\int_{\Omega} \mathcal A\left( (C L_N)^{-1} \sum_{n=1}^N |u_n|^2 \right)dx \leq |\Omega| \,,
$$
where $\mathcal A(t) = e^t - 1 - t$, $L_N=\sum_{n=1}^N n^{-1}$ and where $C$ is a universal constant.

%%%%%%%%%%%%

%%%%%%%%%%%%

\subsection{The Lieb--Thirring inequality}

The original LT inequality in its form for orthonormal functions reads as follows.

\begin{theorem}\label{lt}
	Let $d\geq 1$ and $s>0$. Then, if $u_1,\ldots,u_N\in H^s(\R^d)$ are orthonormal in $L^2(\R^d)$,
	$$
	\sum_{n=1}^N \int_{\R^d} |(-\Delta)^{s/2} u_n|^2\,dx \gtrsim \int_{\R^d} \left( \sum_{n=1}^N \left|u_n \right|^2 \right)^{1+2s/d}dx \,.
	$$
\end{theorem}

\begin{remarks}\label{rem:lt}
	(a) The main point is that the implicit constant can be chosen independently of $N$.\\
	(b) The bound is equivalent to the bound
	$$
	\sum_j |E_j| \lesssim \int_{\R^d} V_-^{1+d/(2s)}dx
	$$
	on the sum of the negative eigenvalues (counted with multiplicities) of the generalized Schr\"odinger operator $(-\Delta)^s+V$ in $L^2(\R^d)$. This equivalence is, for instance, in the sense that the sharp constants in the two inequalities are in one-to-one correspondence.\\
	(c) Theorem \ref{lt} in $d=3$ with $s=1$ is due to Lieb and Thirring \cite{LiTh0} and was the crucial ingredient in their proof of stability of matter; see also \cite{LiSe}. Their proof of Theorem~\ref{liebhls} in \cite{LiTh} for $s=1$ extends to general $s$. Alternative proofs are due to \cite{Ru,LuSo,Sa}.\\
	(d) Lieb and Thirring \cite{LiTh} made a famous conjecture about the optimal constant in the inequality in Theorem \ref{lt} for $s=1$; see, for instance, \cite{Fr} for details. This predicts, in particular, that there is a fundamental difference between dimensions $d\leq 2$ and $d\geq 3$. This conjecture is \emph{open} in any dimension.\\
	(e) The best constants in Theorem \ref{lt} are due to \cite{FrHuJeNa}. A bound with `almost' the semiclassical constant and a gradient remainder term appears in \cite{Nam}.\\
	(f) As a step towards the Lieb--Thirring conjecture, one can study the best constant in the inequality in Theorem \ref{lt} with fixed $N$. For $s=1$ it is shown in \cite{FrGoLe} that in dimensions $d\geq 3$ this constant is always strictly less than the optimal constant that works for arbitrary $N$. This is consistent with the Lieb--Thirring conjecture. For further results in this direction, see also \cite{FrGoLe2,FrGoLe3}. 
\end{remarks}

In the spirit of the generalization discussed in Subsection \ref{sec:hlsgen}, Theorem \ref{lt} has been extended to abstract operators satisfying certain heat kernel bounds or Sobolev inequalities; see \cite{FrLiSe3}.

%%%%%%%%%%%%

\subsection{A more general Lieb--Thirring inequality}

The following theorem provides a Sobolev inequality with exponent $q$ less than $2(1+2s/d)$, the exponent in Theorem~\ref{lt}. The bound is deduced in \cite{LiPa} via a duality argument from a bound of Lieb and Thirring \cite{LiTh}. Note that the functions here are orthogonal, not necessarily orthonormal.

\begin{theorem}\label{ltgen}
	Let $d\geq 1$, $s>0$ and $2<q< 2(1+2s/d)$. Then, if $u_1,\ldots,u_N\in H^s(\R^d)$ are orthogonal in $L^2(\R^d)$,
	$$
	\sum_{n=1}^N \int_{\R^d} |(-\Delta)^{s/2} u_n|^2\,dx \gtrsim \left( \sum_{n=1}^N \|u_n\|_2^{\frac{2(2d-(d-2s)q)}{2d+4s-dq}} \right)^{-\frac{2d+4s-dq}{d(q-2)}} \!\! \left( \int_{\R^d} \!\left( \sum_{n=1}^N \left|u_n \right|^2 \right)^\frac q2 \!dx \right)^\frac{4s}{d(q-2)} \!\!.
	$$
\end{theorem}

\begin{remarks}
	(a) The implicit constant can be chosen independently of $N$.\\
	(b) The bound is equivalent to the bound
	$$
	\sum_j |E_j|^\gamma \lesssim \int_{\R^d} V_-^{\gamma+d/(2s)}dx
	$$
	on the sum of the negative eigenvalues (counted with multiplicities) of the generalized Schr\"odinger operator $(-\Delta)^s+V$ in $L^2(\R^d)$. Here $\gamma>1$ and $q<2(1+2s/d)$ are related by
	$$
	q = \frac{2(\gamma+\frac d{2s})}{\gamma+\frac{d}{2s}-1} \,,
	\qquad
	\gamma = \frac{2d-(d-2s)q}{2s(q-2)} \,.
	$$
	\\
	(c) For $s=1$ Lieb and Thirring \cite{LiTh} made a famous conjecture about the optimal constant in the eigenvalue inequality in (b), which translates into a conjecture for the constant in Theorem \ref{lt}. This conjecture was proved by Laptev and Weidl \cite{LapWei-00} for $\gamma\geq 3/2$, that is, $q\leq 2(d+3)/(d+1)$.\\
	(d) The analysis mentioned in Remark \ref{rem:lt} (f) concerning truncated versions of the inequality is applicable as well in the situation of Theorem \ref{ltgen} with $s=1$; see \cite{FrGoLe,FrGoLe2,FrGoLe3}. 
\end{remarks}

%%%%%%%%%%%%%

\section{Fourier restriction inequalities for orthonormal functions}

We now turn to inequalities for systems of orthonormal functions that are mathematically related to the question of restricting the Fourier transform to hypersurfaces. Such a restriction is possible under certain curvature assumptions on the hypersurface and has important applications to partial differential equations.

\subsection{Strichartz inequality for orthonormal functions}

The Strichartz inequality  \cite{Str,KeTa} concerns solutions $e^{it\Delta}\psi$ of the free Schr\"odinger equation and quantifies their dispersive behavior. It states that if $d\geq 1$, $2\leq q \leq\infty$ if $d=1$, $2\leq q<\infty$ if $d=2$ and $2\leq q \leq 2d/(d-2)$ if $d\geq 3$, and $2/p + d/q=d/2$, then for all $\psi\in L^2(\R^d)$,
\begin{equation}
	\label{eq:strichartz}
	\int_\R \left( \int_{\R^d} |e^{it\Delta}\psi|^q\,dx \right)^{p/q}dt \lesssim \|\psi\|_2^p \,.
\end{equation}
Here is a version of this inequality for systems of orthonormal functions.

\begin{theorem}\label{strichartz}
	Let $d\geq 1$, $2\leq q< 2(d+1)/(d-1)$ and $2/p + d/q = d/2$. Then, if $\psi_1,\ldots,\psi_N\in L^2(\R^d)$ are orthonormal in $L^2(\R^d)$,
	$$
	\int_\R \left( \int_{\R^d} \left( \sum_{n=1}^N \left|e^{it\Delta}\psi_n\right|^2 \right)^{q/2}dx \right)^{p/q}dt \lesssim N^{p(q+2)/(4q)}.
	$$
\end{theorem}

\begin{remarks}\label{strichartzrem}
	(a) The power of $N$ is best possible, as can be deduced from \cite{FrLeLiSe1}.\\
	(b) The bound in Theorem \ref{strichartz} can be slightly improved, namely, for any sequence $0\leq\nu\in\ell^{2q/(q+2)}$,
	$$
	\int_\R \left( \int_{\R^d} \left( \sum_n \nu_n \left|e^{it\Delta}\psi_n\right|^2 \right)^{q/2}dx \right)^{p/q}dt \lesssim \left( \sum_n \nu_n^{2q/(q+2)} \right)^{p(q+2)/(4q)}.
	$$
	The bound in the theorem corresponds to the case $\nu_n\in\{0,1\}$ and is equivalent to a bound for $\nu$ in the Lorentz space $\ell^{2q/(q+2),1}$. It is remarkable that, while in the bounds in the previous section the Lorentz norm $\ell^{s,1}$ on the right side is optimal, here it can be improved to an $\ell^s$ norm. The assumption $\nu\in\ell^{2q/(q+2)}$ cannot be relaxed to $\nu\in\ell^s$ for any $s>2q/(q+2)$ \cite{FrLeLiSe1}.\\
	(c) Theorem \ref{strichartz} appears in \cite{FrLeLiSe1} for $q\leq 2(d+2)/d$ and in \cite{FrSa1} in the full range. The proof in \cite{FrSa1} uses a duality argument, similarly to that in Subsection \ref{sec:duality}. In fact,it is slightly simpler, since the duality between $\mathcal S^r$ and $\mathcal S^{r'}$ is more straightforward than that between $\mathcal S^{r}_{\rm weak}$ and the Lorentz space $\mathcal S^{r',1}$, which is at the core of Subsection \ref{sec:duality}. On the other hand, the fact that here we work in a mixed norm space $L^{p/2}_t L^{q/2}_x$ does not really complicate the argument.\\
	(d) It is conjectured in \cite{FrLeLiSe1} that Theorem \ref{strichartz} remains valid for $q=2(d+1)/(d-1)$. At the same time it is shown there that the strengthening in (b) with the $\ell^{2q/(q+2)}$-norm fails at $q=2(d+1)/(d-1)$. This conjecture was disproved in \cite{BeHoLeNaSa} in dimension $d=1$, but is still \emph{open} for $d\geq 2$.\\
	(e) There is a `semiclassical' version of the inequality where the Schr\"odinger equation is replaced by a transport equation for densities on phase space. The proof in \cite{FrLeLiSe1} can be adapted to this setting, as shown in \cite{BeBeGuLe}. For more on the connection between the two equations, see \cite{Sa0}. The disproof of the conjecture mentioned in (d) for $d=1$ was by disproving the corresponding conjecture in this simpler setting, namely by using the existence of a Kakeya set of zero measure. The validity of this analogue conjecture for $d\geq 2$ is still \emph{open}.\\
	(f) There is a natural `one-particle constant', namely the sharp constant in \eqref{eq:strichartz}. This was determined in the diagonal case $q=p$ in \cite{Fo} for $d=1,2$. Besides, there is a semiclassical constant related to the inequality in (e). To which extent these two constants play a role for the sharp constant in Theorem \ref{strichartz}, in analogy with the Lieb--Thirring conjecture, has not been investigated.
\end{remarks}

The restriction $q<2(d+1)/(d-1)$ in Theorem \ref{strichartz} is not present for the single function inequality \eqref{eq:strichartz}. It is known that for orthonormal functions the case $q\geq 2(d+1)/(d-1)$ behaves differently, but there are several open questions. The following is known.

\begin{theorem}\label{strichartz2}
	Let $d\geq 2$, $2(d+1)/(d-1)\leq q<2d/(d-2)$ and $2/p+d/q=d/2$. Let $\beta< 2q/(d(q-2))$. Then, if $\psi_1,\ldots,\psi_N\in L^2(\R^d)$ are orthonormal in $L^2(\R^d)$,
	$$
	\int_\R \left( \int_{\R^d} \left( \sum_{n=1}^N \left|e^{it\Delta}\psi_n\right|^2 \right)^{q/2}dx \right)^{p/q}dt \lesssim N^{p/(2\beta)}.
	$$
\end{theorem}

\begin{remarks}
	(a) It is known that the bound in the theorem does not hold with an exponent $\beta> 2q/(d(q-2))$, as can be deduced from \cite{FrSa2}, but it is not known whether it holds with exponent $\beta= 2q/(d(q-2))$.\\
	(b) Similarly as in the case of Theorem \ref{strichartz}, the bound in Theorem \ref{strichartz2} can be slightly improved, namely, for any sequence $0\leq\nu\in\ell^\beta$ with $\beta<2q/(d(q-2))$,
	$$
	\int_\R \left( \int_{\R^d} \left( \sum_n \nu_n \left|e^{it\Delta}\psi_n\right|^2 \right)^{q/2}dx \right)^{p/q}dt \lesssim \left( \sum_n \nu_n^\beta \right)^{p/(2\beta)}.
	$$
	The bound is known to fail, in general, if $\nu\in\ell^\beta$ for $\beta>2q/(d(q-2))$ \cite{FrSa2}.\\
	(c) Theorem \ref{strichartz2} appears in \cite{FrSa2} (see the discussion there after Proposition 1). It is obtained by interpolation between Theorem \ref{strichartz} with $q$ near $2(d+1)/(d-1)$ and the bound \eqref{eq:strichartz} with $q$ near its maximal value $2d/(d-2)$.\\
	(d) If the conjecture mentioned in Remark \ref{strichartzrem} (d) is true, an interpolation argument (at least in dimensions $d\geq 3$) might yield Theorem \ref{strichartz2} with $\beta=2q/(d(q-2))$.\\
	(e) Let us discuss the endpoints $q=2d/(d-2)$ in $d\geq 3$ and $q=\infty$ in $d=1$, which are excluded in Theorem \ref{strichartz2}. At the endpoint $q=2d/(d-2)$ in $d\geq 3$, it is known that there is no gain due to orthonormality over the triangle inequality, that is, the bound in (b) holds with $\beta = 1$ and not with any larger power \cite{FrSa2}. At the endpoint $q=\infty$ in $d=1$ it is known that the bound in (b) does not hold for $\beta\geq 2$ (see \cite{FrLeLiSe1} and also \cite[Proposition 1]{FrSa2}) and one may wonder whether it holds for $\beta<2$. In \cite{BeLeNa1} it is shown that the bound holds for $\beta\leq 4/3$ and that the slightly weaker inequality with the $L^{2,\infty}_t L^\infty_x$-norm instead of the $L_t^2 L^\infty_x$-norm holds for all $\beta<2$.
\end{remarks}

Strichartz inequalities for orthonormal system have been proved for more general operators than $-\Delta$ (see, e.g., \cite{FrSa1,BeLeNa2}) and for more regular functions (see, e.g., \cite{BeHoLeNaSa,BeLeNa1,BeLeNa2}).

One application of the Strichartz inequality for orthonormal functions concerns the nonlinear, time-dependent Hartree equation for infinite quantum systems. (Here `infinite' means that the initial data are allowed to have infinite trace.) Using Theorem~\ref{strichartz} one can show global wellposedness and, for small initial data, dispersion for large time; see \cite{FrSa1,Sa0}. For the more involved case of a positive background density, see \cite{LeSa1,LeSa2}.

Another application of the Strichartz inequality for orthonormal functions concerns Besov-space improvements of inequality \eqref{eq:strichartz} for single functions; see \cite[Corollary 9]{FrSa1}. While these bounds can be derived using deep results from bilinear restriction theory, it is interesting to note that the proof via Theorem \ref{strichartz} is much more elementary.

%%%%%%%%%%%%%%%%

\subsection{Stein--Tomas inequality for orthonormal functions}

The Fourier restriction problem is whether the Fourier transform of a function on $\R^d$ has a welldefined restriction to a hypersurface and, if so, to establish corresponding $L^p$ bounds. Sometimes it is helpful to study the equivalent, adjoint problem of Fourier extensions. From a harmonic analysis perspective, the Strichartz inequality corresponds to a Fourier extension inequality for the hypersurface $\{ (\xi,-|\xi|^2) :\ \xi\in \R^d \}$ in $\R^{d+1}$ endowed with a natural measure. Another paradigmatic case concerns the Fourier extension for the sphere. The corresponding result, due to Tomas \cite{To} and Stein \cite{St}, states that, if $f\in L^2(\Sph^{d-1})$, then
\begin{equation}
	\label{eq:steintomas}
	\int_{\R^d} \left| \int_{\Sph^{d-1}} e^{i\omega\cdot x} f(\omega)\,d\omega \right|^{2(d+1)/(d-1)}dx \lesssim \|f\|_{L^2(\Sph^{d-1})}^{2(d+1)/(d-1)} \,.
\end{equation}
The inequality extends trivially to exponents greater than $2(d+1)/(d-1)$ on the left side, but a counterexample due to Knapp shows that $2(d+1)/(d-1)$ is the smallest possible exponent. Here is a version for orthonormal functions.

\begin{theorem}\label{steintomas}
	Let $d\geq 2$. Then, if $f_1,\ldots,f_N$ are orthonormal in $L^2(\Sph^{d-1})$,
	$$
	\int_{\R^d} \left( \sum_{n=1}^N \left| \int_{\Sph^{d-1}} e^{i\omega\cdot x} f_n(\omega)\,d\omega \right|^2 \right)^{(d+1)/(d-1)}dx \lesssim N^{d/(d-1)}.
	$$
\end{theorem}

\begin{remarks}\label{rem:steintomas}
	(a) The power of $N$ is best possible, as can be deduced from \cite{FrSa1}.\\
	(b) Similarly as Theorem \ref{strichartz}, the bound in Theorem \ref{steintomas} can be slightly improved, namely, for any $0\leq\nu\in\ell^{(d+1)/d}$,
	$$
	\int_{\R^d} \left( \sum_n \nu_n \left| \int_{\Sph^{d-1}} e^{i\omega\cdot x} f_n(\omega)\,d\omega \right|^2 \right)^{(d+1)/(d-1)}dx \lesssim \left( \sum_n \nu_n^{(d+1)/d} \right)^{d/(d-1)}.
	$$
	The assumption $\nu\in \ell^{(d+1)/d}$ cannot be relaxed to $\nu\in \ell^r$ for any $r>(d+1)/d$ \cite{FrSa1}.\\
	(c) Theorem \ref{steintomas} appears in \cite{FrSa1}, where it is proved using a duality argument similarly as in Subsection \ref{sec:duality}.\\
	(d) In analogy to Remark \ref{strichartzrem} (f), the optimal constant in \eqref{eq:steintomas} is only known for $d=3$ \cite{Fo2}; see also \cite{FrLiSa} for a connection between the optimal constants in \eqref{eq:steintomas} and \eqref{eq:strichartz}. As far as we know, a `semiclassical inequality' corresponding to that in Theorem \ref{steintomas} has not been investigated.
\end{remarks}

One application of Theorem \ref{steintomas} concerns trace ideal bounds for the scattering matrix for Schr\"odinger operators $-\Delta+V$ in $L^2(\R^d)$ \cite{FrSa1}. These bounds are universal in the sense that they only depend on the `energy' parameter and an $L^p$ norm of $V$, and the trace ideal exponent is shown to be optimal.

To motivate the discussion in the following subsection, we note that by the duality argument in Subsection \ref{sec:duality} and by scaling, the bound in Theorem \ref{steintomas} (or rather in Remark \ref{rem:steintomas} (b)) can be written as
$$
\left\| \mathcal R_k W \right\|_{2(d+1)} \lesssim k^\frac{d-1}{2(d+1)} \left\| W \right\|_{d+1} \,,
$$
where $\mathcal R_k$ denotes restriction of the Fourier transform to the sphere $\{ |\xi|=k\}$. Integrating this bound with respect to $k$ between $\lambda$ and $\lambda+1$, we obtain, in terms of the spectral projection $\Pi_\lambda=\1(\lambda^2\leq-\Delta\leq(\lambda+1)^2)$ with $\lambda\geq 1$,
$$
\left\| \Pi_\lambda |W|^2 \Pi_\lambda \right\|_{d+1} = \left\| \overline W \Pi_\lambda W \right\|_{d+1} \leq \int_\lambda^{\lambda+1} \left\| \overline W \mathcal R_k W \right\|_{d+1}\,dk
\lesssim \lambda^\frac{d-1}{d+1} \left\| W \right\|_{d+1}^2 \,.
$$
Dualizing back, we find that if $(f_n)$ are orthonormal in $L^2(\R^d)$ and satisfy $\supp\widehat f_n \subset\{ \lambda\leq |\xi|\leq\lambda+1\}$ with $\lambda\geq 1$ and if $0\leq\nu\in\ell^{(d+1)/d}$, then
\begin{equation}
	\label{eq:steintomasannulus}
	\left( \int_{\R^d} \left( \sum_n \nu_n |f_n|^2 \right)^\frac{d+1}{d-1}dx \right)^\frac{d-1}{d+1} \lesssim \lambda^\frac{d-1}{d+1} \left( \sum_n \nu_n^\frac{d+1}{d} \right)^\frac{d}{d+1}.
\end{equation}

%%%%%%%%%%%%%%%

\subsection{Spectral cluster bounds}

As shown by Sogge \cite{So}, the version \eqref{eq:steintomasannulus} of the Stein--Tomas inequality has a generalization to closed manifolds. Here is a generalization of this theorem to the case of orthonormal functions from \cite{FrSa3}.

\begin{theorem}\label{clusters}
	Let $(M,g)$ be a smooth compact Riemannian manifold without boundary of dimension $d\geq 2$. Denote by $-\Delta_g$ the Laplace--Beltrami operator on $M$ and, for any $\lambda\ge1$, let $\Pi_\lambda :=\1(\lambda^2\leq-\Delta_g<(\lambda+1)^2)$. Then, if $(f_n)\subset \Pi_\lambda L^2(M)$ are orthonormal in $L^2(M)$ and if $(\nu_n) \subset [0,\infty)$,
	\begin{equation*}
		\left\| \sum_n \nu_n|f_n|^2 \right\|_{L^{q/2}(M)}\lesssim \lambda^{2s(q)}\left(\sum_n \nu_n^{\alpha(q)}\right)^{1/\alpha(q)},
	\end{equation*}
	where
	$$
	\left\{\begin{array}{lll}
		s(q):=d\left(\frac12-\frac1q\right)-\frac12 \,, & \quad \alpha(q)=\frac{q(d-1)}{2d} & \quad \text{if}\ \ \frac{2(d+1)}{d-1}\leq q\leq \infty \,,\\
		s(q):=\frac{d-1}{2}\left(\frac12-\frac1q\right), & \quad \alpha(q)=\frac{2q}{q+2} & \quad \text{if}\ \ 2\leq q\leq\frac{2(d+1)}{d-1} \,.
	\end{array}\right.
	$$
\end{theorem}

\begin{remarks}
	(a) If there is a single nonzero $\nu_n$, the bound in Theorem \ref{clusters} reduces to Sogge's bound \cite{So}. Therefore, according to known results about this inequality, for each $(M,g)$ the exponent $2s(q)$ of $\lambda$ is best possible. As shown in \cite{FrSa3}, for each $(M,g)$ the exponent $\alpha(q)$ is also best possible. Moreover, on $\Sph^2$ with its standard metric it can be shown that the inequality can be saturated even with $\nu_n\in\{0,1\}$ and an arbitrary prescribed sequence $\#\{ n:\ \nu_n =1\}$ \cite{FrSa3}.\\
	(b) The proof of Theorem \ref{clusters} relies on Schatten norm bounds for oscillatory integral operators satisfying the Carleson--Sj\"olin condition, which are of independent interest, but somewhat technical to state. They imply, for instance, Theorem \ref{steintomas}.
\end{remarks}

%%%%%%%%%%%%%%%

\subsection{Kenig--Ruiz--Sogge inequalities}

In this final subsection we discuss resolvent bounds that are close in spirit to the Stein--Tomas theorem. The original result due to Kenig, Ruiz and Sogge \cite{KeRuSo} states that, if $2d/(d+2)\leq p\leq 2(d+1)/(d+3)$ (and $p>1$ if $d=2$), then for all $z\in\C\setminus [0,\infty)$,
\begin{equation}
	\label{eq:kenigruizsogge}
	\left\| (-\Delta-z)^{-1} f \right\|_{p'} \lesssim |z|^{-d/2+d/p-1} \|f\|_p \,.
\end{equation}
For the case $d=2$, see, e.g., \cite{Fr11}. We mention that similar inequalities are valid also on Riemannian manifolds; see, e.g., \cite{DFKeSa,BoShSoYa,FrSc}. A notable feature of the bounds \eqref{eq:kenigruizsogge} is that they do not deteriorate as $z$ approaches the positive real halfline and, for that reason, they are also known as `uniform' Sobolev inequalities. Note that the endpoint exponent $p=2(d+1)/(d+3)$ is the dual of the exponent in the Stein--Tomas Fourier extension inequality \eqref{eq:steintomas} and, in fact, \eqref{eq:steintomas} is an easy consequence of \eqref{eq:kenigruizsogge}.

For $p$ greater than this exponent the uniformity is lost in general. It can be restored, up to $p=2d/(d+1)$ by using mixed norms involving an $L^2$-norm over angular variables \cite{FrSi}. A nonuniform inequality valid for $2(d+1)/(d+3)<p\leq 2$ is
$$
\left\| (-\Delta-z)^{-1} f \right\|_{p'} \lesssim |z|^{-(\frac1p - \frac12)} \dist(z,\R_+)^{-1+(d+1)(\frac1p - \frac12)} \|f\|_p \,.
$$
This bound follows by interpolation between the case $p=2(d+1)/(d+3)$ and the trivial bound at $p=2$. It appeared in an equivalent, dual form in \cite{Fr3}. Remarkably, it is best possible \cite{KwLe}.

Inequality \eqref{eq:kenigruizsogge} is somewhat different from the other ones treated in this paper since it does not involve a Hilbert space norm and, since the operator $(-\Delta-z)^{-1}$ for $z\not\in(-\infty,0]$ is not positive definite, it cannot be rewritten in such a form. Consequently, we cannot state a version for orthonormal functions, but we will directly state trace ideal bounds, similar to what is behind the proofs of the other bounds in this paper. The following two theorems are from \cite{FrSa1} and \cite{Fr3}, respectively.

\begin{theorem}\label{kerusotrace}
	Let $d\geq 2$ and let $8/3\leq q\leq 3$ if $d=2$ and $d\leq q\leq d+1$ if $d\geq 3$. Then, for all $z\in\C\setminus[0,\infty)$,
	$$
	\left\| W_1 (-\Delta-z)^{-1} W_2 \right\|_{(d-1)q/(d-q)} \lesssim |z|^{-1+d/q} \|W_1\|_{q} \|W_2\|_{q} \,.
	$$
\end{theorem}

\begin{theorem}\label{kerusotracenonunif}
	Let $d\geq 1$ and let $q>d+1$. Then, for all $z\in\C\setminus[0,\infty)$,
	$$
	\left\| W_1 (-\Delta-z)^{-1} W_2 \right\|_{q} \lesssim |z|^{-1/q} \dist(z,[0,\infty))^{-1+(d+1)/q} \|W_1\|_{q} \|W_2\|_{q} \,.
	$$
\end{theorem}

The trace ideal exponent in Theorem \ref{kerusotrace} is best possible, as follows from the corresponding result for the Stein--Tomas inequality \cite{FrSa1}. The optimal form of Theorem \ref{kerusotrace} for $d=2$ and $2<q<8/3$ is not known and we refer to \cite{FrSa1} for some partial results.

The main application of Theorems \ref{kerusotrace} and \ref{kerusotracenonunif} is to Lieb--Thirring inequalities for eigenvalues of Schr\"odinger operators with complex-valued potentials; see, e.g., \cite{FrSa1,Fr3}. This is an active area of research with many open question and we refer, for instance, to \cite{FrLaLiSe,Ha,FrLaSa,BoFrVo,BoSt,Cu} for more on this.

%%%%%%%%%%%%%%

\subsection*{Acknowledgements}

Partial support through U.S.~NSF grant DMS-1954995 and through Germany’s Excellence Strategy EXC-2111-390814868 is acknowledged.

%%%%%%%%%%%%%%%%%%%%%%%%%%%%%%%%%%%%%%%%%%%
%%%%%%%%%%%%%%%%%%%%%%%%%%%%%%%%%%%%%%%%%%%

\bibliographystyle{amsalpha}

\end{document}